\documentclass[12pt]{article}

\usepackage{latexsym}
\usepackage{amsmath}
\usepackage{amssymb}
\usepackage{amsthm}

\newtheorem{thm}{Theorem}[section]
\newtheorem{prop}[thm]{Proposition}
\newtheorem{cor}[thm]{Corollary}
\newtheorem{lem}[thm]{Lemma}
\newtheorem{defn}[thm]{Definition}
\newtheorem{rem}[thm]{Remark}

\newtheorem{ex}[thm]{Example}
\newtheorem{nota}[thm]{Notation}

\newcommand{\C}{{\mathbb C}}
\newcommand{\p}{{\mathbb P}}
\newcommand{\g}{{\mathbb G}}
\newcommand{\z}{{\mathbb Z}}

\newcommand{\ideal}{{\cal I}}
\newcommand{\an}{{\cal O}}

\newcommand{\pic}{\operatorname{Pic}}
\newcommand{\lin}{\operatorname{lin}}

\begin{document}

\title{Smooth $n$-dimensional subvarieties of $\p^{2n-1}$ containing a family of very degenerate divisors}
\author{Jos\'e Carlos Sierra \footnote{Partially
supported by the MCYT project BFM2003-03971/MATE and the
Santander-Complutense project PR27/05-13876.}}

\date{\today}
\maketitle

\begin{abstract}
A classification and a detailed geometric description are given for
smooth $n$-dimensional subvarieties $X\subset\p^{2n-1}$ containing a
family of effective divisors each of them spanning a linear $\p^n$
of $\p^{2n-1}$. Some results on multisecant lines to threefolds
$X\subset\p^5$ follow, as a byproduct\-.
\end{abstract}


\smallskip


\section{Introduction}
Let $X\subset\p^r$ be a smooth irreducible $n$-dimensional
projective subvariety defined over the complex numbers. If $r\geq
2n+1$, then $X\subset\p^r$ can be embedded in $\p^{2n+1}$ by a
linear projection, so $X\subset\p^r$ is said to be a low
codimensional subvariety if $r\leq 2n$. Low codimensional
subvarieties are quite special. In the early 70's, Barth and Larsen
showed in \cite{barth}, \cite{bl} and \cite{larsen} that for $i\leq
2n-r$ the cohomology groups $H^i(X,\C)$ and the homotopy groups
$\pi_i(X)$ are the same as those of $\p^r$ (and hence of complete
intersections by Lefschetz hyperplane theorem). According to these
results, Hartshorne \cite{hart} made several conjectures on low
codimensional subvarieties following the principle that the lower
the codimension of $X\subset\p^r$, the closer are its properties to
those of complete intersections. In particular, he conjectured that
$X\subset\p^r$ is a complete intersection if $2r<3n$ and that
$X\subset\p^r$ is linearly normal if $2r-2<3n$. While the first
question remains open, Zak \cite{zak-ln} proved Hartshorne's
conjecture on linear normality by means of the so-called Zak's
Theorem on Tangencies, characterizing also in \cite{zak-sv} the four
Severi varieties as the only non-linearly normal subvarieties
$X\subset\p^r$ in the limit case $2r-2=3n$.

In this paper we study low codimensional subvarieties containing a
family of degenerate divisors by using Zak's Theorem on Tangencies.
Barth-Larsen results imply that the Picard group of $X\subset\p^r$
is generated by its hyperplane section if $r\leq 2n-2$. In this
range, no effective divisor $D$ on $X$ is contained in a linear
$s$-dimensional subspace of $\p^r$ for $n-1\leq s\leq r-2$ (see
Corollary \ref{bl-cor}). So we start by considering the boundary
case $X\subset\p^{2n-1}$, and the simplest situation $D=\p^{n-1}$.

There exist some examples of smooth subvarieties $X\subset\p^{2n-1}$
containing a linear divisor $D=\p^{n-1}$ and realizing any possible
degree allowed by theory (see \cite{bhss}, Ex. 5.2). Moreover the
only smooth subvariety $X\subset\p^{2n-1}$ ruled by a family of
linear divisors is the Segre embedding
$\p^1\times\p^{n-1}\subset\p^{2n-1}$, as Lanteri and Turrini showed
\cite{lt}.

This fact leads us to study possible generalizations. First, one is
tempted to classify smooth subvarieties $X\subset\p^{2n}$ ruled by a
family of linear divisors. A complete classification was given by
Lanteri \cite{lanteri} and Aure \cite{aur} for $n=2$, by Toma
\cite{toma} for $n=3$, and by Ionescu and Toma \cite{it2} for $n\leq
4$. The problem is still open for $n\geq 5$.

Another possibility is to classify smooth subvarieties
$X\subset\p^{2n-1}$ containing a family of effective divisors on $X$
each of them spanning a linear $\p^n$ of $\p^{2n-1}$. This is what
we do in this note. Examples of such varieties were given by Ionescu
and Toma in \cite{it}, where it is shown, for every integer $d\geq
1$, the existence of what they called (degree $d$)-hypersurface
fibrations $X\subset\p^{2n-1}$ over $\p^1$. In this paper we prove
that these are the only hypersurface fibrations in $\p^{2n-1}$, and
we also describe the exceptional cases in which the family of
hypersurfaces does not induce a fibration on $X$. The precise result
is the following.

\bigskip

\noindent{\bf Theorem \ref{main}}\hspace{0.2cm} Let
$X\subset\p^{2n-1}$ be a smooth $n$-dimensional subvariety
containing a family $\{D_t\}_{t\in T}$ of hypersurfaces of degree
$d\geq 2$.
\begin{itemize}
\item [(i)] if $D_t\cap D_{t'}\neq\emptyset$, then $n\leq 3$ with equality if and only if
$X\subset\p^5$ is linked to a linear $\p^3$ by the complete intersection of
a rank-$4$ quadric cone and a hypersurface of $\p^5$ of degree $d$;

\item [(ii)] if $D_t\cap D_{t'}=\emptyset$, then $Y=\bigcup_{t\in T}\langle D_t\rangle\subset\p^{2n-1}$ is a rational normal scroll of
$n$-planes $S_{0,0,1,\ldots,1}\subset\p^{2n-1}$ of degree $n-1$.
Furthermore:
\begin{itemize}
\item[(a)] if $n\geq 3$, then $X\subset\p^{2n-1}$ is linked to a
rational normal scroll of $(n-1)$-planes
$S_{0,0,1,\ldots,1}\subset\p^{2n-3}$ of degree $n-2$ by the complete intersection
of $Y\subset\p^{2n-1}$ and a hypersurface of $\p^{2n-1}$ of degree $d+1$;
\item [(b)] if $n=2$, then $X\subset\p^3$ is a surface of degree $d+1$ containing a
line $L\subset\p^3$ and $\{D_t\}_{t\in T}$ is the pencil $|H-L|$ of plane curves of degree $d$ on $X$.
\end{itemize}
\end{itemize}

\bigskip

To obtain the proof of Theorem \ref{main} we consider two cases
separately. If the intersection of any two divisors of the family is
non-empty, we strongly use Zak's Theorem on Tangencies to get the
bound $n\leq 3$. On the other hand, if the intersection of any two
divisors of the family is empty, we prove that $X\subset\p^{2n-1}$
is linear normal by using Zak's circle of ideas on secant varieties
(see Proposition \ref{ln}).


To complete the overview on related problems, we would like to
stress that smooth $n$-dimensional subvarieties $X\subset\p^{2n}$
containing a family of hypersurfaces are still far from being
classified, even for small $n$. As in the case of subvarieties
containing a family of linear divisors, this is due to the existence
of irregular subvarieties of $\p^{2n}$, in contrast to the case
$r=2n-1$. Consider first $n=2$, i.e. surfaces in $\p^4$ containing a
family of plane curves. There exists a complete classification only
if the plane curves form a fibration on the surface. Ellia and
Sacchiero \cite{es} classified conic fibrations (see also
\cite{br}), and Ranestad \cite{ranestad} classified surfaces ruled
by a family of plane curves of degree greater than or equal to
three. If the plane curves do not form a fibration, just some
partial results have been recently obtained in \cite{st}. As far as
we know, the case $n\geq 3$ remains unexplored. We refer to
\cite{sierra} for a more detailed account on low codimensional
subvarieties containing a family of hypersurfaces.

The last part of the paper is devoted to obtaining some consequences
on multisecant lines to a smooth threefold $X\subset\p^5$ by using
Theorem \ref{main}. Extending results of Severi \cite{severi} and
Aure \cite{aure} on trisecant lines to surfaces in $\p^4$, Kwak
\cite{kwak} proved that a smooth threefold $X\subset\p^5$ has no
apparent triple points if and only if $H^0(\ideal_X(2))\neq 0$. A
classification of threefolds with no apparent quadruple points is
still unknown. We present a different proof of Kwak's result (see
Theorem \ref{triple}) as well as some remarks on the structure of
$4$-secant lines to a threefold $X\subset\p^5$, strengthening both
\cite{mezzetti}, Th. 2.3, and \cite{kwak}, Prop. 4.9. Finally, we
point out that Theorem \ref{main} can also be used to complete the
proof of Th. 1.4 in \cite{mp} (see Remark \ref{gap}).

\section{Preliminaries}

Let $X\subset\p^r$ be an irreducible non-degenerate (i.e. not
contained in a hyperplane) projective subvariety of dimension $n\geq
2$ defined over the complex numbers.

\begin{defn}
\emph{Let $D$ be an effective Weil divisor on $X\subset\p^r$, and
let $\langle D\rangle=\p^s\subset\p^r$ denote its linear span. We
call $D$ a \emph{degenerate divisor on $X$} if $n-1\leq s\leq r-2$.
Moreover $D$ is said to be \emph{linear} if $D=\p^{n-1}$, and $D$ is
said to be a \emph{hypersurface} if $\langle D\rangle=\p^n$. Linear
divisors and hypersurfaces will be called \emph{very degenerate
divisors on $X$}.}
\end{defn}

Let $T$ be an irreducible projective variety. We denote by
$\{D_t\}_{t\in T}$ the algebraic family of divisors on $X$
corresponding to a closed codimension one subvariety ${\cal
D}\subset X\times T$ flat over $T$. Subvarieties $X\subset\p^r$
containing a $2$-dimensional family of very degenerate divisors are
well known. They are described in Lemma \ref{pn} and Theorem
\ref{dim2} below.

\begin{lem}\label{pn}
Let $X\subset\p^r$ be an $n$-dimensional subvariety containing a
family $\{D_t\}_{t\in T}$ of linear divisors on $X$. If $\dim T\geq
2$, then $X=\p^n$.
\end{lem}

\begin{proof}
Since $\dim T\geq 2$, for any two points $x,y\in X$ there exists
some $t\in T$ such that $x,y\in D_t=\p^{n-1}$. Therefore the line
$\langle x,y \rangle$ is contained in $X$, so $X=\p^n$.
\end{proof}

\begin{thm}\label{dim2}
Let $X\subset\p^r$ be an $n$-dimensional subvariety containing a
fa\-mily $\{D_t\}_{t\in T}$ of hypersurfaces. If $r\geq n+2$, then
$\dim T\leq 2$. Moreover if $\dim T=2$, then $\{D_t\}_{t\in T}$ is a
family of (maybe reducible) quadric hypersurfaces and $X\subset\p^r$
is one of the following subvarieties:

\begin{itemize}
\item [(i)] the Segre embedding $\p^1\times\p^2\subset\p^5$;

\item [(ii)] the rational normal scroll $S_{1,2}\subset\p^4$;

\item [(iii)] the Veronese surface $v_2(\p^2)\subset\p^5$ or its projection
to $\p^4$;

\item [(iv)] a cone with vertex $\p^{n-2}$ over a curve;

\item [(v)] a cone over varieties (i), (ii), (iii).
\end{itemize}
\end{thm}

\begin{proof}
If $n=2$, the result was given by C. Segre in \cite{cs}. If $n\geq
3$, we reduce to the case $n=2$ by intersecting $X\subset\p^r$ with
$n-2$ general hyperplanes. If this intersection is a cone with
vertex a point over a curve, then it is easy to see that
$X\subset\p^r$ is a cone with vertex $\p^{n-2}$ over a curve. We
conclude by recalling that the Segre embedding
$\p^1\times\p^2\subset\p^5$ is the only non-trivial extension of
$S_{1,2}\subset\p^4$, and that neither the Veronese surface
$v_2(\p^2)\subset\p^5$ nor its projection to $\p^4$ can be
non-trivially extended. These facts follow, for instance, from
\cite{xxx} and \cite{sd}.
\end{proof}

From now on we also assume that $X\subset\p^r$ is smooth (except in
Section \ref{zak's}). The following result is a corollary of
Barth-Larsen theorems \cite{barth}, \cite{bl} and \cite{larsen}.

\begin{thm}[Barth-Larsen]\label{bl}
Let $X\subset\p^r$ be a smooth subvariety of dimension $n$.

\begin{itemize}
\item [(i)] if $r\leq 2n-1$, then $H^1(X,\an_X)=0$;

\item [(ii)] if $r\leq 2n-2$, then the restriction
$\pic(\p^r)\to\pic(X)$ determines an isomorphism.
\end{itemize}
\end{thm}

\begin{cor}\label{bl-cor}
Let $X\subset\p^r$ be a smooth subvariety of dimension $n$. Let $D$
be a degenerate divisor on $X$. Then $r\geq 2n-1$.
\end{cor}

\begin{proof}
Let $H$ denote the hyperplane divisor on $X\subset\p^r$. Since
$D\subset X$ is a degenerate divisor, there exists an effective
non-zero divisor $E\in |H-D|$ on $X$. If $r\leq 2n-2$, then
$\pic(X)\cong\z$ generated by $H$. Therefore $D\equiv_{\lin} aH$ and
$E\equiv_{\lin} bH$ for some positive integers $a,b$. Hence
$H\equiv_{\lin} D+E\equiv_{\lin} (a+b)H$, yielding a contradiction
(cf. Remark \ref{starting}).
\end{proof}

The bound $r\geq 2n-1$ is sharp, as the following example shows.

\begin{ex}[Beltrametti-Howard-Schneider-Sommese]\label{bhss}
\emph{We present an example of smooth $n$-dimensional subvarieties
$X\subset\p^{2n-1}$ containing a linear $\p^{n-1}$ (see \cite{bhss},
Ex. 5.2). Fix integers $s\geq 0$ and $n\geq 2$. Let ${\cal
P}=\p(\an_{\p^{2n-1}}(1)\oplus\an_{\p^{2n-1}}(s+1))$ and let
$p:{\cal P}\to\p^{2n-1}$ denote the bundle projection. Let $\xi$
denote the tautological line bundle on $\cal P$. The intersection of
$n$ general elements of $|\xi|$ is a smooth $n$-dimensional
subvariety $X'\subset\cal P$ which is mapped isomorphically by $p$
onto its image $X\subset\p^{2n-1}$. The embedding $\cal
P\subset\p^{\dim |\xi|}$ given by $\xi$ contains a linear subspace
$\p^{2n-1}$ correspon\-ding to $\p(\an_{\p^{2n-1}}(1))$, so
$X'\subset\cal P$ (and hence $X\subset\p^{2n-1}$) contains a linear
$\p^{n-1}$. It also follows that $N_{\p^{n-1}/
X}\cong\an_{\p^{n-1}}(-s)$. If $s=0$, then $\deg X=n$ and
$X\subset\p^{2n-1}$ is the Segre embedding
$\p^1\times\p^{n-1}\subset\p^{2n-1}$. If $s\geq 1$, then $\deg
X=\frac{(s+1)^n-1}{s}$. These are the only possible degrees of
smooth $n$-dimensional subvarieties $X\subset\p^{2n-1}$ containing a
linear $\p^{n-1}$, as it is showed in \cite{bhss}, Prop. 8.1.}
\end{ex}

We focus on smooth subvarieties $X\subset\p^{2n-1}$ containing a
family $\{D_t\}_{t\in T}$ of linear divisors.

\begin{rem}\label{vacia}
\emph{Let $X\subset\p^r$ be a smooth $n$-dimensional subvariety
containing a family $\{D_t\}_{t\in T}$ of linear divisors on $X$. If
$D_t\cap D_{t'}\neq\emptyset$ for any $t,t'\in T$, Lemma
\ref{harris} below yields $D_t\cap D_{t'}=\p^{n-2}$. Then it is
straightforward to see that $X=\p^n$ (otherwise $X\subset\p^r$ is a
cone of vertex $\p^{n-2}$, whence singular). Therefore, we can
assume $D_t\cap D_{t'}=\emptyset$.}
\end{rem}

Smooth subvarieties $X\subset\p^{2n-1}$ ruled by a family
$\{D_t\}_{t\in T}$ of linear divisors were first characterized in
\cite{lt}.

\begin{thm}[Lanteri-Turrini]\label{lt}
The only smooth $n$-dimensional subvariety of $\p^{2n-1}$ ruled by a
family of linear divisors is the Segre embedding
$\p^1\times\p^{n-1}\subset\p^{2n-1}$.
\end{thm}

In Section \ref{afterpreliminaries}, we pass on to classify smooth
$n$-dimensional subvarieties $X\subset\p^{2n-1}$ containing a family
$\{D_t\}_{t\in T}$ of hypersurfaces, i.e. effective divisors on $X$
such that $\langle D_t \rangle=\p^n$ for every $t\in T$. First, we
include a few words regarding the main tool we use to obtain our
classification.

\subsection{Zak's Theorem on Tangencies}\label{zak's}

The geometry of projective subvarieties is well reflected in the
relative position of their embedded tangent spaces. Modern algebraic
geometry reduced basically the study of tangents spaces to some
computations involving Chern classes. Fortunately, a new dimension
in the understanding of projective varieties by means of the
relative position of the tangent spaces was introduced by F.L. Zak
and his celebrated theorem on tangencies.

\begin{defn}
\emph{Let $T_xX\subset\p^r$ denote the embedded tangent space to
$X\subset\p^r$ at $x\in X$. A linear subspace $L\subset\p^r$ is said
to be \emph{tangent to $X\subset\p^r$ along $Z\subset X$}, if
$T_xX\subset L$ for every $x\in Z$. We denote $T(Z,X):=\bigcup_{x\in
Z}T_xX\subset\p^r$.}
\end{defn}

\begin{thm}[Zak's Theorem on Tangencies]\label{ztt}
Let $X\subset\p^r$ be a non-degenerate subvariety of dimension $n$.
If a linear subspace $L\subsetneq\p^r$ is tangent to $X\subset\p^r$
along a closed subvariety $Z\subset X$, then
$$\dim Z\leq\dim L-n.$$
\end{thm}

Theorem \ref{ztt} leads to the proof of Hartshorne's conjecture on
linear normality (see \cite{zak} II, Cor. 2.17). Finiteness of Gauss
mappings and the bound on the dimension of the dual variety
(\cite{zak} I, Cor. 2.8 and Cor. 2.5) also follow as an immediate
consequence. References for interesting applications of Theorem
\ref{ztt} by other authors can be found in \cite{zak}. See
\cite{laz} for a brief account on this topic.

\begin{rem}\label{starting}
\emph{Corollary \ref{bl-cor}, which is the starting point of the
paper, can be obtained from Theorem \ref{ztt} as well (see
\cite{zak} I, Prop. 2.16). This second proof does not depend on
Theorem \ref{bl}.}
\end{rem}

\section{Subvarieties $X\subset\p^{2n-1}$ containing a family of
hypersurfaces}\label{afterpreliminaries}

Let $X\subset\p^{2n-1}$ be a smooth $n$-dimensional subvariety
containing a fa\-mily $\{D_t\}_{t\in T}$ of hypersurfaces. Denote
$\p_t^n:=\langle D_t\rangle\subset\p^{2n-1}$ and $Y:=\bigcup_{t\in
T}\p_t^n\subset\p^{2n-1}$.

\begin{rem}\label{segre}
\emph{If $X\subset\p^{2n-1}$ is not the Segre embedding
$\p^1\times\p^2\subset\p^5$, then Theorem \ref{dim2} yields $\dim
T=1$ and $\dim Y=n+1$.}
\end{rem}

\begin{lem}\label{linearly eq}
Let $X\subset\p^{2n-1}$ be a smooth subvariety of dimension $n$. Let
$\{D_t\}_{t\in T}$ be a family of effective divisors on $X$. Then
$D_t$ and $D_{t'}$ are linearly equivalent divisors on $X$ for any
$t,t'\in T$.
\end{lem}

\begin{proof}
It follows from Theorem \ref{bl} that $H^1(X,\an_X)=0$. Then $D_t$
is algebrai\-cally equivalent to $D_{t'}$ if and only if $D_t$ is
linearly equivalent to $D_{t'}$.
\end{proof}

\subsection{Case $D_t\cap D_{t'}\neq\emptyset$}

The next statement, that has already been quoted in Remark
\ref{vacia}, can be found, for instance, in \cite{harris}, Th.
17.24.

\begin{lem}\label{harris}
Let $Z_1,Z_2$ be pure dimensional subvarieties of a smooth variety
$X$. Then the codimension of their intersection (if non-empty) is at
most the sum of their codimensions.
\end{lem}

Lemma \ref{harris} yields $n-2\leq\dim (\p^n_t\cap\p^n_{t'})\leq
n-1$ for general $t,t'\in T$. We analyze both cases separately.

\begin{lem}\label{linear}
Let $X\subset\p^{2n-1}$ be a smooth subvariety of dimension $n$. Let
$\{D_t\}_{t\in T}$ be a family of hypersurfaces such that $D_t\cap
D_{t'}\neq\emptyset$. If $$\dim (\p^n_t\cap\p^n_{t'})=n-1$$ for
general $t,t'\in T$, then $n=2$.
\end{lem}

\begin{proof}
If $\dim (\p^n_t\cap\p^n_{t'})=n-1$ for general $t,t'\in T$, then
either $Y=\p^{n+1}$ (and hence $n=2$) or $Y\subset\p^{2n-1}$ is a
cone with vertex $\p^{n-1}$. In the latter case, it follows from
Remark \ref{segre} that $\dim Y=n+1$. Let $y\in Y$ be a general
point and let $T_yY=\p^{n+1}\subset\p^{2n-1}$ be the tangent space
to $Y$ at $y$. Then $T_yY\subset\p^{2n-1}$ is tangent to $Y$ along
$\langle y,\p^{n-1} \rangle\backslash\p^{n-1}$. Therefore
$T_yY\subset\p^{2n-1}$ is tangent to $X$ along $X_y:=\langle
y,\p^{n-1} \rangle\cap X$. Since $\dim X_y=n-1$, it follows from
Theorem \ref{ztt} that $n\leq 2$.
\end{proof}

\begin{lem}\label{n3}
Let $X\subset\p^{2n-1}$ be a smooth subvariety of dimension $n$. Let
$\{D_t\}_{t\in T}$ be a family of hypersurfaces such that $D_t\cap
D_{t'}\neq\emptyset$. If $$\dim (\p^n_t\cap\p^n_{t'})=n-2$$ for
general $t,t'\in T$, then $n=3$.
\end{lem}

\begin{proof}
Fix general $t,t'\in T$. Since $\dim (\p^n_t\cap\p^n_{t'})=n-2$,
then $D_t\cap D_{t'}=\p^{n-2}$. If $X\subset\p^{2n-1}$ is not the
Segre embedding $\p^1\times\p^2\subset\p^5$, the family
$\{D_t\}_{t\in T}$ of hypersurfaces is the pencil of divisors on $X$
generated by $D_t$ and $D_{t'}$. This follows from Remark
\ref{segre} and Lemma \ref{linearly eq}. In particular,
$\p^{n-2}\subset D_t$ for every $t\in T$ and $Y\subset\p^{2n-1}$ is
an $(n+1)$-dimensional cone with vertex $\p^{n-2}$. As before, the
tangent space $T_yY=\p^{n+1}\subset\p^{2n-1}$ to $Y$ at a general
$y\in Y$ is tangent to $X$ along $X_y:=\langle y,\p^{n-2}\rangle\cap
X$. Since $\dim X_y=n-2$, Theorem \ref{ztt} yields $n\leq 3$.
Therefore $n=3$.
\end{proof}

\begin{nota}
\emph{Let $a_0\leq\ldots\leq a_n$ be a non-trivial sequence of
non-negative integers. We denote by
$S_{a_0,\ldots,a_n}\subset\p^{n+a_0+\cdots +a_n}$ the rational
normal scroll of $n$-planes of degree $a_0+\cdots +a_n$. For
$E=\oplus_{i=0}^n\an_{\p^1}(a_i)$, we can also see
$S_{a_0,\ldots,a_n}\subset\p^{n+a_0+\cdots +a_n}$ as the image of
$\p(E)$ in $\p^{n+a_0+\cdots +a_n}$ by the tautological line bundle
$\an_{\p(E)}(1)$ on $\p(E)$. Recall that the rational normal scroll
of $n$-planes $S_{1,\ldots,1}\subset\p^{2n+1}$ of degree $n+1$ is
precisely the Segre embedding $\p^1\times\p^n\subset\p^{2n+1}$.}
\end{nota}

\begin{prop}\label{quadric}
Let $X\subset\p^{2n-1}$ be a smooth $n$-dimensional subvariety
containing a family $\{D_t\}_{t\in T}$ of hypersurfaces such that
$D_t\cap D_{t'}\neq\emptyset$. Then $n\leq 3$, with equality if and
only if $X\subset\p^5$ is linked to a linear $\p^3$ by the complete
intersection of a rank-4 quadric cone and a hypersurface of $\p^5$.
\end{prop}

\begin{proof}
Lemmas \ref{linear} and \ref{n3} yield $n\leq 3$. The Segre
embedding $\p^1\times\p^2\subset\p^5$ is linked to a linear $\p^3$
by the complete intersection of a rank-4 quadric cone and a quadric
hypersurface of $\p^5$. If $X\subset\p^5$ is not this Segre
embedding, it follows from Remark \ref{segre} and the proof of Lemma
\ref{n3} that $Y\subset\p^5$ is a rational scroll of $3$-planes of
vertex $L=\p^1$. Since $X\subset\p^5$ is linearly normal by
\cite{zak-ln}, and $D_t\cap D_{t'}=L$, it follows that
$Y\subset\p^5$ is a rational normal scroll of $3$-planes. Hence
$\deg Y=2$ and $Y\subset\p^5$ is the quadric cone
$S_{0,0,1,1}\subset\p^5$. The result follows, for instance, from
\cite{kw}, Th. 2.1.
\end{proof}

\subsection{Case $D_t\cap D_{t'}=\emptyset$}
In this case, the main ingredient of Theorem \ref{main} is the
following proposition.

\begin{prop}\label{ln}
Let $X\subset\p^{2n-1}$ be a smooth subvariety of dimension $n$. Let
$\{D_t\}_{t\in T}$ be a family of hypersurfaces. If $D_t\cap
D_{t'}=\emptyset$, then $X\subset\p^{2n-1}$ is linearly normal.
\end{prop}

\begin{proof}
Suppose to the contrary that $X\subset\p^r$ is a non-degenerate
$n$-dimen\-sional subvariety that can be isomorphically projected to
$\p^{2n-1}$. We denote by $SX\subset\p^r$ the secant variety of
$X\subset\p^r$. Let $z\in SX$ be a general point and let $E_z\subset
X$ be the corresponding entry locus, i.e. the set of points in $X$
swept out by secant lines to $X\subset\p^r$ passing through $z$.
Since $\dim SX=2n-1$, then $\dim E_z=2$. Terracini's Lemma \cite{tl}
yields $T_xX\subset T_zSX$ for every $x\in E_z$. Let
$J(D_t,D_{t'})\subset\p^r$ denote the join of $D_t$ and $D_{t'}$,
i.e. the set of points contained in lines spanned by a point of
$D_t$ and a point of $D_{t'}$. We claim that $J(D_t,D_{t'})$ has the
expected dimension $$\dim J(D_t,D_{t'})=2(n-1)+1=2n-1.$$ Otherwise,
projecting to $\p^{2n-2}$ from a suitable linear subspace of $\p^r$,
we get $\emptyset = D_t\cap D_{t'}\subset\p^{2n-2}$, which is not
possible. For every $t\in T$, denote $Z_t=E_z\cap D_t$. Since $\dim
J(D_t,D_{t'})=\dim SX$, it follows that $J(D_t,D_{t'})=SX$ for any
$t,t'\in T$, $t\neq t'$. Therefore, for every $t\in T$, there exist
infinitely many secant lines to $X$ passing through $z$ and meeting
$D_t$. Thus $\dim Z_t\geq 1$ for every $t\in T$, and so $\dim Z_t=1$
for a general $t\in T$, since $\dim E_z=2$. Theorem \ref{ztt} yields
$T(Z_t,D_t)=\p^n_t\subset\p^{2n-1}$, and hence
$$Y=\bigcup_{t\in T}\p^n_t=\bigcup_{t\in
T}T(Z_t,D_t)\subset\bigcup_{t\in T}T(Z_t,X)=T(E_z,X)\subset T_zSX,$$

\noindent so $Y\subset T_zSX$. But $T_zSX\subset\p^r$ is a linear
subspace of dimension $2n-1$. Hence $X\subset\p^r$ is a degenerate
subvariety, a contradiction.
\end{proof}

\begin{lem}\label{y}
Let $X\subset\p^{2n-1}$ be a smooth subvariety of dimension $n$. Let
$\{D_t\}_{t\in T}$ be a family of hypersurfaces. If $D_t\cap
D_{t'}=\emptyset$, then $Y\subset\p^{2n-1}$ is a rational normal
scroll of $n$-planes $S_{0,0,1,...,1}\subset\p^{2n-1}$ of degree
$n-1$.
\end{lem}

\begin{proof}
$X\subset\p^{2n-1}$ is a linearly normal subvariety by Proposition
\ref{ln}. It follows from Theorem \ref{dim2} and Lemma \ref{linearly
eq} that $T$ is a rational curve. Since $D_t\cap D_{t'}=\emptyset$,
we deduce that $Y\subset\p^{2n-1}$ is a rational normal scroll of
$n$-planes. Therefore $\deg Y=n-1$. We claim that
$\dim(\p^n_t\cap\p^n_{t'})=1$ for any $t,t'\in T$. Otherwise, if
$\dim(\p^n_t\cap\p^n_{t'})\geq 2$, then $D_t\cap
D_{t'}\neq\emptyset$. Thus we obtain
$Y=S_{0,0,1,...,1}\subset\p^{2n-1}$.
\end{proof}

\begin{ex}[Ionescu-Toma]
\emph{For every integer $d\geq 1$, the existence of smooth
$n$-dimensional (degree $d$)-hypersurface fibrations
$X\subset\p^{2n-1}$ contained in $Y=S_{0,0,1,...,1}\subset\p^{2n-1}$
was shown in \cite{it}, Prop. 3. We reproduce their construction
(some notations have been changed). Let $E=\an_{\p^1}^{\oplus
2}\oplus\an_{\p^1}(1)^{\oplus n-1}$ and let $\p(E)$ be the
projective bundle over $\p^1$. Let $\pi:\p(E)\to\p^1$ be the
projection, $H=\an_{\p(E)}(1)$ and $F=\pi^*\an_{\p^1}(1)$. It
follows from Bertini's theorem that the general element $X\in|dH+F|$
is smooth. Let $\varphi:\p(E)\to Y\subset\p^{2n-1}$ be the morphism
induced by $H$. Then the restriction $\varphi:X\to\p^{2n-1}$ turns
out to be an embedding onto its image. We prove in the sequel that
moreover elements $X\in|dH+F|$ are the only smooth divisors on
$\p(E)$ such that $\varphi:X\to\p^{2n-1}$ is an embedding.}
\end{ex}

\begin{prop}\label{key} With the above notations, let $X\in|dH+bF|$ be a smooth divisor on
$\p(E)$. If $\varphi:X\to\p^{2n-1}$ is an embedding, $X\notin |F|$, then $b=1$.
\end{prop}

\begin{proof}
Set $E'=\an_{\p^1}^{\oplus 2}$, $E''=\an_{\p^1}(1)^{\oplus n-1}$ and
$S=\p(E')$. Let $p_1=\pi_{|S}$ and $p_2=\varphi_{|S}$ be the two
projections of $S\cong\p^1\times\p^1$. Since $S$ is mapped by
$\varphi$ onto the vertex $L=\p^1$ of $Y$, it follows that
$\varphi:\p(E)\backslash S\to Y\backslash L\subset\p^{2n-1}$ is an
embedding. Therefore $\varphi:X\to\p^{2n-1}$ is an embedding only if
the restriction $\varphi:X\cap S\to L$ is an embedding. Then
$C=X\cap S$ is a rational curve on $S$ and $C\in |dF_1+bF_2|$, where
$F_i=p_i^*(\an_{\p^1}(1))$. Since $\varphi$ contracts $F_1$, we get
$1=F_1\cdot C=b$.
\end{proof}

\begin{rem}\label{flavor}
\emph{To get a more projective flavor, we present a different
geome\-tric description of $X\in|dH+F|$ embedded in
$Y\subset\p^{2n-1}$. Let $R\in |H-F|$. Since $H\equiv_{\lin} R+F$,
it follows that $dH+F\equiv_{\lin}(d+1)H-R$. If $n\geq 3$, then
$\varphi(R)\subset Y$ is a rational normal scroll of $(n-1)$-planes
$S_{0,0,1,\ldots,1}\subset\p^{2n-3}$ of degree $n-2$. If $n=2$, then
$\varphi(R)=L\subset\p^3$. Rational normal scrolls are projectively
normal (see \cite{ah}, Prop. 2.9), i.e. the restriction maps
$$H^0(\an_{\p^{2n-1}}(m))\to H^0(\an_{Y}(m))$$ are surjective for
every integer $m\geq 1$. Therefore if $n\geq 3$, then
$X\subset\p^{2n-1}$ is linked to a rational normal scroll of
$(n-1)$-planes $S_{0,0,1,\ldots,1}\subset\p^{2n-3}$ of degree $n-2$
by the complete intersection of $Y$ and a hypersurface of
$\p^{2n-1}$ of degree $d+1$. On the other hand, if $n=2$, then
$X\subset\p^3$ is a surface of degree $d+1$ containing $L=\p^1$, and
$\{D_t\}_{t\in T}$ is the pencil $|H-L|$ of plane curves of degree
$d$ on $X$.}
\end{rem}

\begin{rem}\label{n=3}
\emph{Let us look more closely at $n=3$. In this case
$Y=S_{0,0,1,1}\subset\p^5$ is a rank-$4$ quadric cone. Hence
$Y\subset\p^5$ is ruled by two pencils of linear $\p^3$'s, say $T_1$
and $T_2$. If, for instance, $X\subset\p^5$ is linked to a $\p^3$ of
$T_2$ by the complete intersection of $Y$ and a hypersurface of
$\p^5$ of degree $d+1$, then $\{D_t\}_{t\in T_1}$ induces a
fibration on $X$ by hypersurfaces of degree $d$. On the other hand,
$\{D_t\}_{t\in T_2}$ is a pencil of hypersurfaces of degree $d+1$
with base locus the vertex $L=\p^1$ of $Y\subset\p^5$, as in
Proposition \ref{quadric}.}
\end{rem}

The main result of the paper is now obtained.

\begin{thm}\label{main}
Let $X\subset\p^{2n-1}$ be a smooth $n$-dimensional subvariety
containing a family $\{D_t\}_{t\in T}$ of hypersurfaces of degree
$d\geq 2$.
\begin{itemize}
\item [(i)] if $D_t\cap D_{t'}\neq\emptyset$, then $n\leq 3$ with equality if and only if
$X\subset\p^5$ is linked to a linear $\p^3$ by the complete intersection of
a rank-$4$ quadric cone and a hypersurface of $\p^5$ of degree $d$;

\item [(ii)] if $D_t\cap D_{t'}=\emptyset$, then $Y=\bigcup_{t\in T}\langle D_t\rangle\subset\p^{2n-1}$ is a rational normal scroll of
$n$-planes $S_{0,0,1,\ldots,1}\subset\p^{2n-1}$ of degree $n-1$.
Furthermore:
\begin{itemize}
\item[(a)] if $n\geq 3$, then $X\subset\p^{2n-1}$ is linked to a
rational normal scroll of $(n-1)$-planes
$S_{0,0,1,\ldots,1}\subset\p^{2n-3}$ of degree $n-2$ by the complete intersection
of $Y\subset\p^{2n-1}$ and a hypersurface of $\p^{2n-1}$ of degree $d+1$;
\item [(b)] if $n=2$, then $X\subset\p^3$ is a surface of degree $d+1$ containing a
line $L\subset\p^3$ and $\{D_t\}_{t\in T}$ is the pencil $|H-L|$ of plane curves of degree $d$ on $X$.
\end{itemize}
\end{itemize}
\end{thm}

\begin{proof}
If follows from Propositions \ref{quadric}, \ref{ln}, \ref{key},
Lemma \ref{y} and Remark \ref{flavor}.
\end{proof}

For $n\geq 3$, the Segre embedding
$\p^1\times\p^{n-1}\subset\p^{2n-1}$ is linked to a rational normal
scroll of $(n-1)$-planes $S_{0,0,1,\ldots,1}\subset\p^{2n-3}$ of
degree $n-2$ by the complete intersection of $Y\subset\p^{2n-1}$ and
a quadric hypersurface of $\p^{2n-1}$. Thus we can state Theorems
\ref{lt} and \ref{main} together if $n\geq 3$.

\begin{cor}\label{together}
Let $X\subset\p^{2n-1}$ be a smooth $n$-dimensional subvariety
containing a family of very degenerate divisors on $X$. If $n\geq
3$, then $X\subset\p^{2n-1}$ is linked to a rational normal scroll
of $(n-1)$-planes $S_{0,0,1,\ldots,1}\subset\p^{2n-3}$ of degree
$n-2$ by the complete intersection of a rational normal scroll of
$n$-planes $S_{0,0,1,\ldots,1}\subset\p^{2n-1}$ of degree $n-1$ and
a hypersurface of $\p^{2n-1}$.
\end{cor}

\begin{rem}\emph{Proposition \ref{key} can be obtained as a particular case of \cite{ilic}, Th. 3.7.
In fact, using the terminology of \cite{ilic}, we can reformulate
Corollary \ref{together} as follows: a smooth $n$-dimensional
subvariety $X\subset\p^{2n-1}$ contains a family of very degenerate
divisors if and only if $X\subset\p^{2n-1}$ is a linearly normal
\emph{Roth variety}. We wish to thank Paltin Ionescu for pointing
out this paper.}
\end{rem}

\begin{rem}
\emph{In particular, smooth $n$-dimensional subvarieties
$X\subset\p^{2n-1}$ containing a family of hypersurfaces are
linearly normal. This is no longer true if either $X\subset\p^{2n}$
or $\langle D_t \rangle=\p^{n+1}$, as the following example shows.}
\end{rem}

\begin{ex}
\emph{The Veronese surface $v_2(\p^2)\subset\p^5$ contains a family
$\{D_t\}_{t\in T}$ of conics and can be isomorphically projected to
$\p^4$. The Segre embedding $\p^2\times\p^2\subset\p^8$ contains a
family $\{D_t\}_{t\in T}$ of divisors $\p^1\times\p^2\subset\p^5$
and can be isomorphically projected to $\p^7$. In both cases
$D_t\cap D_{t'}\neq\emptyset$. On the other hand, if $D_t\cap
D_{t'}=\emptyset$, it is easy to see that subvarieties
$X\subset\p^{2n}$ ruled by a family of linear divisors are linearly
normal. We do not know if hypersurface fibrations $X\subset\p^{2n}$
are linearly normal as well. At least, this is known to be true for
$n\leq 3$ (cf. \cite{fujita}).}
\end{ex}

\section{Multisecant lines to threefolds in $\p^5$}

Let $\Sigma_k(X)\subset\g(1,r)$ be the variety of $k$-secant lines
to a smooth irreducible non-degenerate projective subvariety
$X\subset\p^r$, and let $S_k(X)\subset\p^r$ be the closure of the
union of all $k$-secant lines. If $X\subset\p^5$ is a threefold,
then $S_2(X)=\p^5$. Moreover, it follows from the
(dimension+2)-secant lemma \cite{ran} that $S_5(X)\neq\p^5$. We say
that $X\subset\p^5$ has no apparent triple (resp. quadruple) points
if $S_3(X)\neq\p^5$ (resp $S_4(X)\neq\p^5$). Threefolds in $\p^5$
with no apparent triple points were classified by Kwak \cite{kwak},
extending to $n=3$ the following result of Aure \cite{aure}.

\begin{thm}[Aure]\label{aaaure}
Let $X\subset\p^4$ be a smooth surface. Then $S_3(X)\neq\p^4$ if
and only if either $H^0(\ideal_X(2))\neq 0$ or $X\subset\p^4$ is a
quintic elliptic scroll.
\end{thm}

We reobtain Kwak's result in Theorem \ref{triple}, as a byproduct of
Theorem \ref{main} and the following lemma (see \cite{bs}).

\begin{lem}[B. Segre]\label{bs}
Let $V\subset\p^r$ be an irreducible subvariety of dimension $n$.
Let $\Sigma_{\infty}\subset\g(1,r)$ be a component of maximal
dimension of the variety of lines contained in $V$. Then:
\begin{enumerate}
\item [(i)] if $\dim\Sigma_{\infty}=2n-2$, then $V=\p^n$;

\item [(ii)] if $\dim\Sigma_{\infty}=2n-3$, then $V\subset\p^r$ is either a quadric hypersurface or
a scroll of linear $\p^{n-1}$'s over a curve.
\end{enumerate}
\end{lem}

\begin{thm}[Kwak]\label{triple}
Let $X\subset\p^5$ be a smooth threefold. Then $S_3(X)\neq\p^5$ if
and only if $H^0(\ideal_X(2))\neq 0$.
\end{thm}

\begin{proof}
If $X\subset\p^5$ is contained in a quadric hypersurface
$Q\subset\p^5$, then $S_3(X)\subset Q$. For the converse, let
$C\subset\p^3$ be the intersection of $X\subset\p^5$ with a general
$\p^3$ of $\p^5$. If $\dim\Sigma_3(X)\leq 4$, then $C\subset\p^3$ is
a smooth curve with at most a finite number of trisecant lines. So
$C\subset\p^3$ is either a rational normal cubic or the complete
intersection of two quadrics of $\p^3$. Therefore $X\subset\p^5$ is
either the Segre embedding $\p^1\times\p^2\subset\p^5$ or the
complete intersection of two quadric hypersurfaces of $\p^5$. In
both cases $H^0(\ideal_X(2))\neq 0$. So we assume
$\dim\Sigma_3(X)=5$. Then $S_3(X)\subset\p^5$ is either a quadric
hypersurface or a scroll of $\p^3$'s over a curve by Lemma \ref{bs}.
If the second occurs, then $S_3(X)\subset\p^5$ is a quadric cone
$S_{0,0,1,1}\subset\p^5$ by Theorem \ref{main}. So $X\subset\p^5$ is
contained in the quadric hypersurface $S_3(X)\subset\p^5$.
\end{proof}

\begin{cor}
Let $X\subset\p^{n+2}$ be a smooth subvariety of dimension $n$. Then
$S_3(X)\neq\p^{n+2}$ if and only if either $H^0(\ideal_X(2))\neq 0$
or $X\subset\p^4$ is a quintic elliptic scroll.
\end{cor}

\begin{proof}
Let $C\subset\p^3$ denote the intersection of $X\subset\p^{n+2}$
with $n-1$ general hyperplanes. If $\dim\Sigma_3(X)\leq 2n-2$, then
$C\subset\p^3$ is a smooth curve with at most a finite number of
trisecant lines and $X\subset\p^{n+2}$ is either a complete
intersection of two quadric hypersurfaces or $X\subset\p^5$ is the
Segre embedding $\p^1\times\p^2\subset\p^5$. If
$\dim\Sigma_3(X)=2n-1$, then $S_3(X)\subset\p^{n+2}$ is either a
quadric hypersurface or a scroll of $\p^n$'s over a curve by Lemma
\ref{bs}. In the scroll case, $X\subset\p^{n+2}$ contains a family
of degenerate divisors, so Corollary \ref{bl-cor} yields $n\leq 3$.
Thus the result follows from Theorems \ref{aaaure} and \ref{triple}.
\end{proof}

A classification of threefolds $X\subset\p^5$ with no apparent
quadruple points is still an open problem. If
$H^0(\ideal_X(3))\neq 0$, then $X\subset\p^5$ has no apparent
quadruple points but it is not known if there exists such a
threefold satisfying $H^0(\ideal_X(3))=0$ (see \cite{kwak} and
\cite{mezzetti}). Theorem \ref{main} is useful to sharpen
\cite{mezzetti}, Th. 2.3.

\begin{prop}
Let $X\subset\p^5$ be a smooth threefold. If $\dim\Sigma_k(X)=5$ for
some integer $k\geq 4$, then $H^0(\ideal_{X}(2))\neq 0$.
\end{prop}

\begin{proof}
We claim that $S_k(X)\neq\p^5$. To the contrary, assume
$S_k(X)=\p^5$. Let $H\subset\p^4$ be a general hyperplane. Denote
$S=X\cap H\subset\p^4$. Let $p\in H$ be a general point. Since
$\dim\Sigma_k(X)=5$, there exists a $1$-dimensional family of
$k$-secant lines to $X\subset\p^5$ passing through $p$. Thus $p\in
H$ is contained in some $k$-secant line to $S\subset\p^4$,
contradicting \cite{ran}. Therefore $S_k(X)\subset\p^5$ is either a
quadric hypersurface or a scroll of $\p^3$'s over a curve. Hence the
result follows from Theorem \ref{main}.
\end{proof}

\begin{rem}\label{gap}
\emph{It was stated in \cite{mp}, Th. 1.4, that a smooth threefold
$X\subset\p^5$ containing a family $\{D_t\}_{t\in T}$ of surfaces of
$\p^3$ is contained in a quadric hypersurface of $\p^5$. In the
proof, the case $D_t\cap D_{t'}\neq\emptyset$ was not considered.
This fact was communicated by Emilia Mezzetti to Angelo Felice
Lopez, who pointed it out to the author. Thus Theorem \ref{main}
also helps to fill this gap.}
\end{rem}

Finally, we slightly improve Prop. 4.9 in \cite{kwak} by means of
Theorem \ref{main}.

\begin{prop}\label{multi}
Let $X\subset\p^5$ be a smooth threefold. If $S_4(X)\neq\p^5$ and
$H^0(\ideal_X(3))=0$, then $\dim\Sigma_4(X)=4$.
\end{prop}

\begin{proof}
If $\dim\Sigma_4(X)=5$, then $S_4(X)\subset\p^5$ is either a quadric
hypersurface or a scroll of $\p^3$'s over a curve by Lemma \ref{bs}.
Hence $S_4(X)\subset\p^5$ is a quadric hypersurface by Theorem
\ref{main}, so $H^0(\ideal_X(3))\neq 0$. Now assume
$\dim\Sigma_4(X)\leq 3$. Then $\deg X\leq 8$ by the proof of Prop.
4.9 in \cite{kwak}. Smooth threefolds $X\subset\p^5$ up to degree
$8$ are described in \cite{ion}, \cite{ion2}, \cite{okonek} and
\cite{okonek2}. They are contained in a cubic hypersurface of $\p^5$
except Palatini scroll, which has one apparent quadruple point.
Therefore $\dim\Sigma_4(X)=4$.
\end{proof}

\begin{rem}
\emph{It was stated in \cite{kwak}, Prop. 4.9, that
$S_4(X)\subset\p^5$ is a hypersurface if $H^0(\ideal_X(3))=0$ and
$\deg X\geq 9$, but no proof is given for this assertion when
$\dim\Sigma_4(X)=4$. If $\dim\Sigma_4(X)=4$ and $\dim S_4(X)\leq
3$, then each of the $4$-dimensional irreducible components of
$\Sigma_4(X)\subset\g(1,5)$ sweeps out a linear $\p^3$ of $\p^5$
by Lemma \ref{bs}. This might happen if $X\subset\p^5$ contains a
finite number of surfaces of degree at least four, each of them
contained in some $\p^3$.}
\end{rem}

Jos\'e Carlos Sierra

Departamento de \'Algebra

Facultad de Ciencias Matem\'aticas

Universidad Complutense de Madrid

28040 Madrid, Spain

e-mail: jcsierra@mat.ucm.es
\end{document}